\documentclass[12pt, a4paper]{amsart}
\usepackage{tikz}
\usepackage{amsmath,amsthm,amsfonts,amssymb,latexsym,mathrsfs,color,extarrows}
\usepackage{color, graphicx, comment}

\allowdisplaybreaks

\newtheorem{Theorem}{Theorem}
\newtheorem{Corollary}[Theorem]{Corollary}
\newtheorem{Proposition}[Theorem]{Proposition}

\newtheorem{Lemma}[Theorem]{Lemma}
\newtheorem{Question}[Theorem]{Question}
\theoremstyle{definition}

\font\KFracFont=cmr12 at 20pt
\def\KKK{\mathop{\lower 4pt\hbox{\KFracFont K}}\limits}

\newcommand{\Sym}{\mathfrak{S}}

\DeclareMathOperator{\des}{des}

\DeclareMathOperator{\maj}{maj}

\DeclareMathOperator{\inv}{inv}
\DeclareMathOperator{\coinv}{coinv}

\DeclareMathOperator{\sumt}{sum}
\DeclareMathOperator{\len}{len}
\DeclareMathOperator{\id}{id}

\usepackage{tikz}

\usepackage{epsf}

\usepackage[english]{babel}

\title{Cyclotomic Euler-Mahonian polynomials}

\author{Guo-Niu HAN}
\address{I.R.M.A., UMR 7501, Universit\'e de Strasbourg
et CNRS, 7 rue Ren\'e Descartes, F-67084 Strasbourg, France}
\email{guoniu.han@unistra.fr}



\date{September 16, 2026}

\begin{document}

\begin{abstract} 
The cyclotomic Eulerian polynomials and the cyclotomic Mahonian polynomials 
have each been the subject of extensive studies in Combinatorics, 
with particular attention to their signed versions.
In contrast, the joint study of cyclotomic Euler-Mahonian polynomials 
has received far less consideration.
To the best of our knowledge, the only prior result in this direction 
is a formula due to Wachs for the signed Euler-Mahonian polynomials 
in the even case.
In this paper, we focus on the cyclotomic Euler-Mahonian polynomials
and derive a formula based on the Hadamard product.
	As corollaries, we obtain the $I$-analogue (where $I=\sqrt{-1}$) of Wachs' formula 
for signed Euler-Mahonian polynomials, as well as the previously missing 
odd case for the signed Euler-Mahonian polynomials.

\end{abstract}

\maketitle

\section{Introduction}  
Let $A_n(t,q,p)$ denote the generating function for the symmetric group $\Sym_n$ with respect to the trivariate statistic $(\des,\maj,\inv)$, namely
$$
A_n(t,q,p) = \sum_{\sigma\in \mathfrak{S}_n} t^{\des(\sigma)} q^{\maj(\sigma)} p^{\inv(\sigma)}.
$$
The specialization $p=1$, written $A_n(t,q):=A_n(t,q,1)$, recovers the usual {\it Euler–Mahonian polynomials}. In addition, the polynomials $A_n(t,1)$ and $A_n(1,q)$ are known as the {\it Eulerian polynomials} and the {\it Mahonian polynomials}, respectively.

\medskip

For every positive integer $m$, let $\xi_m$ denote a primitive $m$-th root of unity. By analogy with the case $m=2$ (so that $\xi_2=-1$), which gives rise to the {\it signed} polynomials, we call the specialization $p=\xi_m$ the {\it cyclotomic} polynomials of order $m$. The cyclotomic Eulerian and cyclotomic Mahonian polynomials have been widely studied individually in combinatorics, with particular emphasis on the signed Eulerian and signed Mahonian 
cases \cite{Adin2005GR_sigMahonians, Desarmenien1983_funcsym, Desarmenien1992F_qEuler, Loday1989_op,  Wachs1992_involution}. A more detailed historical overview of this line of work is provided in Section \ref{sec:notation}.

\medskip

However, the cyclotomic Euler-Mahonian polynomials have received little attention so far. 
To the best of our knowledge, the only result in this area is due to Wachs \cite{Wachs1992_involution}, who derived a formula for the signed Euler-Mahonian polynomials in the even case, $A_{2n}(t,q,-1)$; see Theorem \ref{th:Wachs}.

In this paper, we investigate the joint cyclotomic Euler-Mahonian polynomials. Let $(t;q)_n := (1-t)(1-tq)\cdots(1-tq^{n-1})$ denote the standard $q$-ascending factorial. Recall that the {\it Hadamard product} of two power series is given by
$$
\sum_{\ell \geq 0} a_\ell t^{\ell}
\otimes
\sum_{\ell \geq 0} b_\ell t^{\ell}
=
\sum_{\ell \geq 0} a_\ell b_\ell t^{\ell}.
$$
By employing the Hadamard product, we derive an explicit expression for the cyclotomic Euler-Mahonian polynomials. The main theorem is presented below.
\begin{Theorem}\label{th:main}
Let $\xi_m$ denote a primitive $m$-th root of unity. For any positive integer $n$,
express $n$ in the form $n = mk + i$ with $0 \leq i \leq m - 1$. Then
\begin{equation}\label{eq:main}
\frac{A_{n}(t,q,\xi_m)}{(t;q)_{n+1}}
=
\frac{A_i(t,q,\xi_m)}{(t;q)_{i+1}}
\otimes
\frac{A_k(t,q^m)}{(t;q^m)_{k+1}}.
\end{equation}
\end{Theorem}

\medskip

Theorem \ref{th:main} naturally specializes to several previously known identities for 
the polynomials $A_n(t,q,\xi_m)$ by setting 
specific values for the parameters  $t,q,m$ and $i$. 
In Section \ref{sec:cases}, we recover the  following identities: 

\begin{itemize}

\item Case $q=1$ and $m=2$: This recovers the results for the signed Eulerian polynomials,
originally studied by
D\'esarm\'enien,  Foata, Loday and Wachs 
\cite{ Desarmenien1983_funcsym,  Desarmenien1992F_qEuler,  Loday1989_op,  Wachs1992_involution}.
See Theorem \ref{th:DF}.

\item	Case $t=1$ and $m=2$: 
This yields the identities for the signed Mahonian polynomials,
originally studied
by Gessel, Simion and Wachs \cite{Wachs1992_involution}. 
See Theorem \ref{th:GesselSimion}.


\item case $t=1$: 
This recovers the identity for the cyclotomic Mahonian polynomials,
studied by Adin, Gessel  and Roichman
		\cite{Adin2005GR_sigMahonians}. See Theorem \ref{th:AdinGR}.

\item case $m=2, i=0$: 
This provides the even case of
the signed Euler-Mahonian polynomials, studied by Wachs 
		\cite{Wachs1992_involution}. See Theorem \ref{th:Wachs}.

\end{itemize}We also derived several further corollaries, which are presented in Section \ref{sec:cases}. 
We highlight just two of them: one is an $I$-analogue (where $I=\sqrt{-1}$) of Wachs’s formula for signed Euler–Mahonian polynomials, and the other provides the previously missing odd case for signed Euler–Mahonian polynomials.
\begin{Corollary}\label{cor:caseI}
Let $I=\sqrt{-1}$, we have 
\begin{equation}\label{eq:caseI}A_{4k}(t,q,I) 
=	 \frac{  (t;q)_{4k+1}}{(t;q^4)_{k+1}} A_k(t,q^4). 
\end{equation}
\end{Corollary}

\begin{Corollary}\label{cor:missing} We have
\begin{align}
A_{2k+1}(t,q,-1) 
&=(t;q)_{2k+2} \sum_{\ell\geq 0}  [\ell+1]_q \, [\ell+1]_{q^2}^k 
					 t^\ell.\label{eq:missing}
\end{align}

\end{Corollary}

Section \ref{sec:notation} reviews the basic preliminary concepts and the original results concerning the signed Eulerian and Mahonian polynomials. We present two proofs of the main theorem. The first follows the approach developed by Adin, Gessel, and Roichman \cite{Adin2005GR_sigMahonians}, see Section \ref{sec:proof1}. The second proof is more direct, see Section \ref{sec:proof2}. Both proofs make use of Gessel’s generating function formula for $A_n(t,q,p)$ \cite{Gessel1977}.

\section{Preliminaries and Notation}\label{sec:notation} 
Let $\Sym_n$ denote the symmetric group on the set $\{1,2,\ldots,n\}$. For a permutation 
$\sigma = \sigma_1 \sigma_2 \cdots \sigma_n \in \Sym_n$,
the {\it inversion number}, {\it descent number}, and {\it major index} are defined respectively by
\begin{align*}
\inv \sigma &= \#\{(i,j)\mid 1 \le i < j \le n,\ \sigma_i > \sigma_j\},\\
\des \sigma &= \#\{i \mid 1 \le i \le n-1,\ \sigma_i > \sigma_{i+1}\},\\
\maj \sigma &= \sum \{i \mid 1 \le i \le n-1,\ \sigma_i > \sigma_{i+1}\},
\end{align*}
where the notation $``\#"$ denotes the cardinality of a set.
Recall that $A_n(t,q,p)$ denotes the generating polynomial for the symmetric group $\Sym_n$ with respect to the statistics $(\des, \maj, \inv)$,
$$
A_n(t,q,p) = \sum_{\sigma\in \mathfrak{S}_n} t^{\des(\sigma)} q^{\maj(\sigma)} p^{\inv(\sigma)}.
$$
We abbreviate $A_n(t,q):=A_n(t,q,1)$ for the Euler–Mahonian polynomials. The $q$-{\it ascending factorial} is defined by
\begin{align*}
	(a;q)_n &= \begin{cases}
	1, &\text{if $n=0$};\\
	(1-a)(1-aq)\cdots(1-aq^{n-1}), &\text{if $n\ge 1$};
\end{cases}\\
	(a;q)_\infty &= \prod_{n\ge 0} (1-aq^n).
\end{align*}
For brevity, we shall also use throughout the paper the shorthand notation
$$
(q)_n := (q;q)_n,
$$
and likewise for any other variable, so that, for instance, $(p)_n=(p;p)_n$.
Write $[n]_q=1+q+\cdots + q^{n-1}$ and $[n]_q! = [1]_q [2]_q \cdots [n]_q$. 
We define the $q$-binomial coefficient by
$$
{n \brack k}_q = \frac{(q;q)_n}{(q;q)_k (q;q)_{n-k}}.
$$
It is a classical result \cite{Andrews1976, MacMahon1913} that $A_n(1,q) = [n]_q!$ and that
\begin{equation}\label{eq:qbinom}
\frac{1}{(t;q)_n} = \sum_{k\geq 0} {n+k-1 \brack k}_q\, t^k
\quad\text{and}\quad
\frac{1}{(t;q)_\infty} = \sum_{k\geq 0} \frac{t^k}{(q)_k}.
\end{equation}
Carlitz \cite{Carlitz1954, Carlitz1975} obtained the following identity for the Euler–Mahonian polynomials.

\begin{Theorem}[Carlitz]\label{th:Carlitz}
For all $n\geq 0$,
\begin{equation}\label{eq:Carlitz}
\frac{A_n(t,q)}{(t;q)_{n+1}} = \sum_{\ell \ge 0} t^\ell \big([\ell+1]_q\big)^n.
\end{equation}
\end{Theorem}


D\'esarm\'enien and Foata established the following identities linking the Eulerian polynomials with the signed Eulerian polynomials, confirming a conjecture of Loday \cite{ Desarmenien1983_funcsym,  Desarmenien1992F_qEuler,  Loday1989_op}.
\begin{Theorem}[Loday, D\'esarm\'enien, Foata]\label{th:DF}
We have
\begin{align}
A_{2k}(t,1,-1) &={  (1-t)^{k}} A_k(t,1),\label{eq:DFeven}\\
A_{2k+1}(t,1,-1) &=(1-t)^{k} {A_{k+1}(t,1) }.\label{eq:DFodd}
\end{align}
\end{Theorem}

Wachs  \cite{ Wachs1992_involution} provided  a bijective proof of Theorem \ref{th:DF}. 
In addition, her bijection leads to a formula of Gessel and Simion for the signed Mahonian polynomials.
\begin{Theorem}[Gessel, Simion, Wachs]\label{th:GesselSimion}
We have
\begin{equation}\label{eq:GesselSimion}A_{n}(1,q, -1) =\left(\frac{1-q  }{ 1+q}\right)^{\lfloor \frac{n}{2}\rfloor} [n]_q!.
\end{equation}
\end{Theorem}

Moreover, Wachs's bijection gives a $q$-analogue of \eqref{eq:DFeven}, corresponding to the even case of Theorem \ref{th:DF}, for the signed Euler–Mahonian polynomials.
\begin{Theorem}[Wachs]\label{th:Wachs} We have
\begin{equation}\label{eq:Wachs}A_{2k}(t,q,-1) =	(tq;q^2)_k   A_k(t,q^2).
\end{equation}
\end{Theorem}
As Wachs observed, her bijection does not yield a $q$-analogue in the odd case for signed Euler–Mahonian polynomials.

\medskip

The cyclotomic Mahonian polynomials were studied by Adin, Gessel, and Roichman \cite{Adin2005GR_sigMahonians}. In their work, they derived a generalization of Theorem \ref{th:GesselSimion} of Gessel and Simion.

\begin{Theorem}[Adin, Gessel, Roichman]\label{th:AdinGR}
Let $n = mk + i$ with $0 \leq i \leq m - 1$. Then
\begin{equation}\label{eq:AdinGR}
A_n(1, q, \xi_m) = \frac{(q;q)_n}{(q;q)_i (1-q^m)^k} A_i(1, q, \xi_m).
\end{equation}
\end{Theorem}

The cyclotomic Eulerian polynomials have also been analyzed by D\'esarm\'enien and Foata, who expressed $A_n(t,1,\xi_m)$ in terms of a two-index Eulerian polynomial; see \cite{Desarmenien1983_funcsym, Desarmenien1992F_qEuler}.
\medskip

Since then, these signed enumeration identities have been broadened and developed along two principal lines. On the one hand, they have been extended from the symmetric group to other reflection groups by Adin–Gessel–Roichman \cite{Adin2005GR_sigMahonians}, 
Reiner \cite{Reiner1995_Des},  Biagioli \cite{Biagioli2006_sMah}, Biagioli–Caselli \cite{Biagioli2004C}, Brenti–Sentinelli \cite{Brenti2020S}, Eu–Fu–Lo \cite{Eu2022FL}, and Eu–Lin–Lo \cite{Eu2021LL_sEulMah}. On the other hand, the focus has shifted from general permutations to various families of restricted permutations, including 321-avoiding permutations studied by Adin–Roichman \cite{Adin2004R} and Eu–Fu–Pan–Ting \cite{Eu2015FPT}, simsun permutations examined by Eu–Fu–Pan \cite{Eu2014FP_simsun}, PRW permutations studied by Lin–Wang–Zeng \cite{Lin2019WZ}, permutations with subsequence constraints considered by Eu–Fu–Hsu–Liao–Sun \cite{Eu2020FHLS}, and derangements analyzed by Ji–Zhang \cite{Ji2025Z_sDerang}.
Finally, let us mention the recent work of Dong and Lin \cite[Theorem 1.11]{Dong2024L}, which
provides a refinement of D\'esarm\'enien--Foata's formula \eqref{eq:DFeven}--\eqref{eq:DFodd} by
the peak statistic.

\medskip

\section{Specializations}\label{sec:cases} 
In this section, we specialize our main theorem by assigning specific values to the parameters
$t,q,m$ and $i$.
In doing so, we recover several identities already stated in Section \ref{sec:notation},
and we further obtain additional special cases presented as corollaries.
We begin by rewriting \eqref{eq:main} in a different form.
For each $i$ with $0\le i\le m-1$, let $C_i(\ell; q, \xi_m)$ $(\ell\ge 0)$ be the coefficients
defined by the expansion
\begin{equation}\label{eq:case_small_i}\frac{A_i(t,q,\xi_m)}{ (t;q)_{i+1}} =
\sum_{\ell \geq 0} C_i(\ell; q, \xi_m) t^\ell. 
\end{equation}
Using Carlitz’s formula \eqref{eq:Carlitz}, the main theorem is equivalent to 
\begin{equation}\label{eq:main2}\frac{A_{n}(t,q,\xi_m)}{ (t;q)_{n+1}} 
= 
\sum_{\ell \geq 0} C_i(\ell; q, \xi_m) \, ([\ell+1]_{q^m})^k t^\ell,
\end{equation}
providing $n=mk+i$ with $0\le i\le m-1$.

\subsection{Case $i=0$} 
In this case we have $n=mk$.
Since $A_0(t,p,q)=1$, the main theorem gives
\begin{equation*}
\frac{A_{n}(t,q,\xi_m)}{ (t;q)_{n+1}} 
= 
	\frac{1}{1-t}
\otimes
\frac{A_k(t,q^m)}{(t;q^m)_{k+1}} 
= 
\frac{A_k(t,q^m)}{(t;q^m)_{k+1}}.
\end{equation*}
Thus, in this situation we obtain an explicit expression without resorting to the Hadamard product.
\begin{Corollary}\label{cor:case_i0} We have
\begin{equation}\label{eq:case_i0}A_{mk}(t,q,\xi_m) 
	=	 \frac{  (t;q)_{mk+1}}{(t;q^m)_{k+1}} A_k(t,q^m). 
\end{equation}
\end{Corollary}

Taking $m=2$, i.e., $\xi_2=-1$ in \eqref{eq:case_i0}, we recover
Wachs' identity stated in Theorem \ref{th:Wachs} for the
signed Euler--Mahonian polynomials. Similarly, setting $m=4$, namely $\xi_4 = I=\sqrt{-1}$ in \eqref{eq:case_i0}, we obtain Corollary \ref{cor:caseI}. 

\subsection{Case $i=1$} 
Here we have $n=mk+1$.
Since $A_1(t,q,p)=1$, we get
\begin{equation*}
\frac{A_i(t,q,\xi_m)}{ (t;q)_{i+1}} 
=\frac{1}{ (t;q)_{2}} 
= \sum_{\ell\geq 0}[\ell+1]_q  t^\ell.
\end{equation*}
Applying the main theorem, we obtain:
\begin{Corollary}\label{cor:case_i1} We have
\begin{align}
A_{mk+1}(t,q,\xi_m) 
&=(t;q)_{mk+2} \sum_{\ell\geq 0}  [\ell+1]_q  \, [\ell+1]_{q^m}^k
					 t^\ell.\label{eq:case_i1}
\end{align}
\end{Corollary}
For $m=2$, this yields Corollary \ref{cor:missing}, which supplies the missing odd case for the signed Euler–Mahonian polynomials.

\subsection{Case $i=2$} 
In this situation we assume $m \ge 3$ and $n = mk+2$.  Since $A_2(t,q,p)=1+tpq$, we obtain
\begin{align*}
\frac{A_2(t,q,\xi_m)}{(t;q)_{3}}
&=\frac{1+tq\xi_m}{(t;q)_3}\\
&=(1+tq\xi_m)\sum_{\ell\ge 0} {2+\ell \brack \ell}_q\,t^\ell\\
&=\sum_{\ell\ge 0} \frac{(1-q^{\ell+2})+q\xi_m(1-q^{\ell})}{1-q^2}\,
[\ell+1]_q\,t^\ell.
\end{align*}
Applying the main theorem \eqref{eq:main2}, we obtain:

\begin{Corollary}\label{cor:case_i2}
For $n=mk+2$ we have
\begin{align*}
\frac{A_{mk+2}(t,q,\xi_m)}{(t;q)_{mk+3}}
&=\sum_{\ell\ge 0}\frac{(1-q^{\ell+2})+q\xi_m(1-q^{\ell})}{1-q^2}\,
[\ell+1]_q\, [\ell+1]_{q^m}^k\,t^\ell.
\end{align*}
\end{Corollary}

\subsection{Case $t=1$}
This specializes to the Mahonian polynomial case.
Recall first the following useful result, the “$t=1$ Lemma” \cite[p. 163]{Foata2004H}.

\begin{Lemma}\label{lem:t=1}
Let $(a_\ell)_{\ell\ge 0}$ be a sequence of formal series with limit $a$, in the sense that
$o(a-a_\ell)\to\infty$ as $\ell\to\infty$. Define
\[
b(t):=(1-t)\sum_{\ell\ge 0} a_\ell\,t^\ell.
\]
Then $b(1)=a$.
\end{Lemma}

We now apply the “$t=1$ Lemma” to identity \eqref{eq:main2}.
Multiplying both sides by $1-t$ gives
\begin{align*}
\frac{A_n(t,q,\xi_m)(1-t)}{(t;q)_{n+1}}
&=(1-t)\sum_{\ell\ge 0}C_i(\ell;q,\xi_m)
   \frac{(1-q^{(\ell+1)m})^k}{(1-q^m)^k}\,t^\ell.
\end{align*}
Setting $t=1$ and invoking Lemma \ref{lem:t=1}, we obtain
\begin{equation}\label{eq:An1}
\frac{A_n(1,q,\xi_m)}{(q;q)_n}
= \frac{1}{(1-q^m)^k}\lim_{\ell\to\infty}C_i(\ell;q,\xi_m).
\end{equation}
Now choose $n=mk+i$ with $k=0$, that is, $n=i$. Then
\begin{equation}\label{eq:limB}
\frac{A_i(1,q,\xi_m)}{(q;q)_i}
= \lim_{\ell\to\infty}C_i(\ell;q,\xi_m).
\end{equation}
Substituting \eqref{eq:limB} into \eqref{eq:An1} yields
\begin{equation*}
A_n(1,q,\xi_m)
= \frac{(q;q)_n}{(q;q)_i(1-q^m)^k}\,A_i(1,q,\xi_m).
\end{equation*}
This is precisely Theorem \ref{th:AdinGR} of Adin, Gessel, and Roichman.
Furthermore, when $m=2$, we recover  Theorem \ref{th:GesselSimion} due to Gessel–Simion.
\subsection{Case $q=1$} 
This corresponds to the Eulerian polynomial setting.
By substituting $q=1$ into Corollaries \ref{cor:case_i0}, \ref{cor:case_i1} and \ref{cor:case_i2}, we obtain:

\begin{Corollary}\label{cor:q=1} We have
\begin{align*}
A_{mk}(t,1,\xi_m) &=(1-t)^{(m-1)k} A_k(t,1),\\
A_{mk+1}(t,1,\xi_m) &=(1-t)^{(m-1)k} A_{k+1}(t,1),\\
A_{mk+2}(t,1,\xi_m) &=(1-t)^{(m-1)k}\Big(\tfrac{1+\xi_m}{2}\,A_{k+2}(t,1)\\
&\qquad\qquad\qquad
+\tfrac{1-\xi_m}{2}\,(1-t)\,A_{k+1}(t,1)\Big).
\end{align*}
\end{Corollary}

\begin{proof}
The first equality follows immediately from setting $q=1$ in \eqref{eq:case_i0}.
For the second, substitute $q=1$ into \eqref{eq:case_i1}; then, using Carlitz’s formula \eqref{eq:Carlitz}, we obtain
\begin{align*}
A_{mk+1}(t,1,\xi_m) 
&=(1-t)^{mk+2} \sum_{\ell\geq 0} (\ell+1)^{k+1} t^\ell\\
&=(1-t)^{mk+2} \frac{A_{k+1}(t,1)}{(1-t)^{k+2}}.
\end{align*}
For the third identity, substituting $q=1$ in Corollary \ref{cor:case_i2} gives
\begin{align*}
A_{mk+2}(t,1,\xi_m) 
&=(1-t)^{mk+3}\sum_{\ell\geq 0} \frac{(\ell+2) + \ell\,\xi_m}{2}\,
(\ell+1)^{k+1} t^\ell\\
&=(1-t)^{mk+3}\sum_{\ell\geq 0} \left(\tfrac{1+\xi_m}{2}(\ell+1)^{k+2}
+\tfrac{1-\xi_m}{2}(\ell+1)^{k+1}\right) t^\ell\\
&=(1-t)^{mk+3}\left(\frac{1+\xi_m}{2}\frac{A_{k+2}(t,1)}{(1-t)^{k+3}}\right.\\
&\qquad\qquad\qquad\left.
+\frac{1-\xi_m}{2}\frac{A_{k+1}(t,1)}{(1-t)^{k+2}}\right).\qedhere
\end{align*}
\end{proof}

In the special case $m=2$, this reduces to the theorem of Loday, D\'esarm\'enien and Foata, Theorem~\ref{th:DF}.

\subsection{Case $m=4$} 
In this situation, $\xi_m = I = \sqrt{-1}$. Substituting $m = 4$ and $\xi_m = I$ into the previously discussed cases $t = 1$ and $q = 1$ respectively, we obtain 
\begin{align*}
A_{4k}(t,1,I) &={  (1-t)^{3k}} A_k(t,1),\\
A_{4k+1}(t,1,I) &=(1-t)^{3k} {A_{k+1}(t,1) },\\
A_{4k+2}(t,1,I) &=(1-t)^{3k}\Big(\tfrac{1+I}{2}A_{k+2}(t,1)\\
&\qquad\qquad\quad +\tfrac{1-I}{2}(1-t)A_{k+1}(t,1)\Big);
\end{align*}
and
\begin{align*}
A_{4n}(1,q, I) &=\frac{ (1-q)^{4n} }{ (1-q^4)^n} [4n]_q!,\\
 {A_{4n+1}(1,q, I)} 
&=\frac{(1-q)^{4n+1}}{(1-q)(1-q^4)^n } [4n+1]_q!,\\
 {A_{4n+2}(1,q, I)} 
&=\frac{(1+qI)}{(1-q)(1-q^2)} \frac{(1-q)^{4n+2}}{(1-q^4)^n } [4n+2]_q!,\\
 {A_{4n+3}(1,q, I)} 
&=\frac{ ( 1 -q - q^2)   +  I q(1 +q -q^2)    }{(1-q) (1-q^2)(1-q^3)} \frac{(1-q)^{4n+3}}{(1-q^4)^n } [4n+3]_q!.
\end{align*}


\section{First proof of the main theorem}\label{sec:proof1} 
Recall the following
generating function for the trivariate statistic $(\des, \maj, \inv)$ obtained by Gessel \cite{Gessel1977}:

\begin{Theorem}[Gessel]\label{th:Gessel}
We have 
\begin{align}
\sum_{n\geq 0}  \frac{ A_n(t,q,p) u^n}{(p;p)_n  (t;q)_{n+1}} 
	&=\sum_{\ell\geq 0} \frac{t^\ell}{\prod_{i=0}^{\ell}\prod_{j\geq 0} (1- u q^{i} p^{j}) }.\label{eq:Gessel}
\end{align}
\end{Theorem}

We now adopt the approach developed by Adin, Gessel, and Roichman \cite{Adin2005GR_sigMahonians}, phrased in terms of cyclotomic polynomials.  
Consider the sequence of polynomials $f_0(p), f_1(p), \dots$ in $p$, defined through the \textit{Eulerian generating function}
\begin{equation}
F(u;p) = \sum_{n=0}^{\infty} f_n(p) \frac{u^n}{(p;p)_n}. \label{eq:F}
\end{equation}
Let $\Phi_m(p)$ denote the \textit{cyclotomic polynomial} of order $m$ in the variable $p$.  
For polynomials $f(p)$ and $g(p)$ with rational coefficients, we write $f(p) \equiv g(p) \pmod{\Phi_m(p)}$ if and only if $f(\xi_m) = g(\xi_m)$.
Let 
$$
F(u;p) = \sum_{n=0}^\infty f_n(p)\frac{u^n}{(p;p)_n} 
\text{\ and\ } 
G(u;p) = \sum_{n=0}^\infty g_n(p)\frac{u^n}{(p;p)_n}
$$ 
be two Eulerian generating functions.  
We then define the congruence $F(u;p) \equiv G(u;p)$ to mean that the corresponding coefficients satisfy $f_n(p) \equiv g_n(p) \pmod{\Phi_m(p)}$ for all $n \geq 0$. 
The key properties of these congruences are summarized in the following lemma, due to
Adin, Gessel, and Roichman \cite{Adin2005GR_sigMahonians}:
\begin{Lemma} \label{lem:AdinGR}
Let $v_i := u^i/(p)_i$.
\begin{enumerate}
    \item If $0 \leq i, j < m$ and $i+j \geq m$, then $v_i v_j \equiv 0$.
    \item If $0 \leq i < m$, then $ v_{mk+i} \equiv {v_{m}^k} v_i / {k!} $.
	\item We have $(p)_{m-1} \equiv m$.
\end{enumerate}
\end{Lemma}
Now we are prepared to prove the main theorem, making use of Theorem \ref{th:Gessel}.

\begin{proof}[Proof of Theorem \ref{th:main}] Define
\begin{equation*}
L(\ell)
= \frac{1}{  \prod_{i=0}^{\ell}\prod_{j\geq 0}   (1- u q^{i} p^{j}) },
\end{equation*}
so that
\begin{equation}\label{eq:AL}
\sum_{n\geq 0}  \frac{ A_n(t,q,p) v_n}{ (t;q)_{n+1}} 
	= \sum_{\ell\geq 0} L(\ell) t^\ell.  
\end{equation}
We compute
\begin{align*}
L(\ell)
&=\prod_{i=0}^{\ell}\prod_{j\geq 0}  \exp \left(\log \frac{1}{   (1-uq^{i} p^{j})  } \right)
	= \prod_{d\geq 1} \exp  \left(\sum_{i=0}^{\ell}\sum_{j\geq 0}  \frac{(uq^ip^j)^d}{d}\right)\\
&=\prod_{d\geq 1} \exp \left( \frac{ (1-q^{(\ell+1)d})  u^d}{d(1-q^d)(1-p^d)}\right)\\
&\equiv \prod_{d= 1}^{m-1} \exp(L'(d))\times  \exp(L'(m)) \times \prod_{d\geq m+1} \exp(L'(d)),
\end{align*}
where
$$
	L'(d):= \frac{ (1-q^{(\ell+1)d})  (p)_{d-1}  }{d(1-q^d) } v_d.
$$
We now examine these three factors separately. For the third factor, we immediately obtain
$$
\prod_{d\geq m+1} \exp(L'(d))  \equiv 1.
$$For the second factor, using Lemma \ref{lem:AdinGR} (2) and (3), we have
\begin{align*}
\exp(L'(m))
&=\exp\left( \frac{(1-q^{(\ell+1)m})(p)_{m-1}v_m}{m(1-q^m)} \right)
\equiv \exp\left( \frac{(1-q^{(\ell+1)m})v_m}{1-q^m} \right) \\
&=\sum_{k\geq 0} \frac{(1-q^{(\ell+1)m})^k v_m^k}{(1-q^m)^k k!}
\equiv \sum_{k\geq 0} ([\ell+1]_{q^m})^k v_{mk}.
\end{align*}
By Lemma \ref{lem:AdinGR} (1), the first factor is a polynomial in $u$. Thus we may write
\begin{equation*}
\prod_{d=1}^{m-1} \exp(L'(d))
\equiv \sum_{i=0}^{m-1} C_i(\ell; q,p)\, v_i.
\end{equation*}
Combining all three factors, we obtain
\begin{align*}
L(\ell)
&\equiv
\left( \sum_{i=0}^{m-1} C_i(\ell; q,p)\, v_i \right)
\left( \sum_{k\geq 0} ([\ell+1]_{q^m})^k v_{mk} \right).
\end{align*}
Therefore, using identity \eqref{eq:AL}, we get
\begin{align*}
\sum_{n\geq 0} \frac{A_n(t,q,p)\, v_n}{(t;q)_{n+1}}
&\equiv \sum_{\ell\geq 0}
\left( \sum_{i=0}^{m-1} C_i(\ell; q,p)\, v_i
\sum_{k\geq 0} ([\ell+1]_{q^m})^k v_{mk} \right) t^\ell.
\end{align*}
Now, comparing coefficients of $v_n$ on both sides for $n=mk+i$ and using Lemma \ref{lem:AdinGR} (2), we obtain $v_{mk} v_i / v_{mk+i} = 1$, and hence
\begin{align}
\frac{A_n(t,q,p)}{(t;q)_{n+1}}
&\equiv \sum_{\ell\geq 0} C_i(\ell; q,p)\,
([\ell+1]_{q^m})^k\, t^\ell.
\label{eq:ABl}
\end{align}
First, substituting $k=0$ and $n=i$ for $0\le i\le m-1$ into \eqref{eq:ABl}, we obtain
\begin{equation*}
\sum_{\ell\geq 0}  C_i(\ell; q,p)\, t^\ell
\equiv \frac{A_i(t,q,p)}{(t;q)_{i+1}}.
\end{equation*}
Next, applying Carlitz’s formula \eqref{eq:Carlitz}, we have
\[
 \sum_{\ell\geq 0} ([\ell+1]_{q^m})^k\, t^\ell 
 = \frac{A_k(t,q^m)}{(t;q^m)_{k+1}}.
\]
Therefore, identity \eqref{eq:ABl} can be rewritten as
\begin{equation*}
\frac{A_{n}(t,q,p)}{(t;q)_{n+1}} 
\equiv
\frac{A_i(t,q,p)}{(t;q)_{i+1}} 
\otimes
\frac{A_k(t,q^m)}{(t;q^m)_{k+1}}.
\end{equation*}
This completes the proof of the main theorem.
\end{proof}
\section{Second proof of the main theorem}\label{sec:proof2} 
Throughout this section, the summation condition
$$
\sum_{n_0, \ldots, n_\ell \vdash n}
$$
will be a shorthand notation for a sum ranging over all nonnegative integers 
$n_j$ $(0\leq j \leq \ell)$ 
satisfying $n_0+n_1 + \cdots + n_\ell =n$.
Recall that the $q$-multinomial coefficient is defined by
$${n_0+n_1+\cdots + n_\ell \brack n_0, n_1, \ldots, n_\ell}_q = \frac{(q;q)_{n_0+n_1+\cdots + n_\ell} }{(q;q)_{n_0} (q;q)_{n_1}\cdots (q;q)_{n_\ell}}.$$
Starting from Gessel's identity \eqref{eq:Gessel} and using the second expansion
in \eqref{eq:qbinom}, with $t$ and $q$ replaced by $uq^i$ and $p$, we obtain
\begin{align*}
[t^\ell] \sum_{n\geq 0}  \frac{ A_n(t,q,p) u^n}{(p;p)_n  (t;q)_{n+1}} 
	&=\prod_{i=0}^{\ell} \frac{1}{\prod_{j\geq 0} (1- (u q^{i}) p^{j}) }\\
&=\prod_{i=0}^{\ell} \sum_{n_i \geq 0} \frac{(uq^i)^{n_i}}{(p)_{n_i}}\\
&=\sum_{n\geq 0} \sum_{n_0, \ldots, n_{\ell}\vdash n}q^{n_1 + 2 n_2 + \cdots + \ell n_{\ell} }  {n   \brack n_0, \ldots, n_{\ell}}_p \frac{u^{n}}{(p;p)_n}.
\end{align*}
By comparing the coefficient of $v_n=u^n/(p)_n$, we deduce the following result. This formulation has the benefit that it can be directly specialized to $p=\xi_m$, in contrast to the original version of Gessel's theorem.
\begin{Theorem}\label{th:util} We have
\begin{align}
[t^{\ell}] \frac{A_n(t,q,p)}{(t;q)_{n+1}} 
	&=\sum_{n_0,  \cdots, n_{\ell} \vdash n}  q^{n_1 + 2 n_2 + \cdots + \ell n_{\ell} }  {n   \brack n_0, \cdots, n_{\ell}}_p.\label{eq:util}
\end{align}
\end{Theorem}

We now present a self-contained proof of \eqref{eq:util}, based on MacMahon’s method \cite{MacMahon1913}.
Fix a positive integer $n$. Let $\mathcal{W}_n=\{0,1,2,\ldots \}^n$ denote the set of all words over the nonnegative integers of length $n$, and let $\mathcal{P}_n$ be the set of all partitions $\lambda = (\lambda_1 \geq \lambda_2 \geq \cdots)$ with largest part $\lambda_1\leq n$.
Let $\mathcal{D}_n$ denote the set of all weakly decreasing words in $\mathcal{W}_n$, and,
for $\overline w\in \mathcal{D}_n$, let $\Sym(\overline w)$ denote the set of all
rearrangements of $\overline w$, so that $\mathcal{W}_n$ is the disjoint union of the
sets $\Sym(\overline w)$ with $\overline w\in \mathcal{D}_n$.
For any word (or permutation, or partition) $w= x_1 x_2 \cdots x_n$, define
$\max(w)=\max\{x_1, x_2, \ldots, x_n\}$ and
$\sumt(w)=x_1 + x_2 + \cdots + x_n$.
If $w$ is a word, define 
\begin{align*}
\inv w&=\#\{(i,j) \mid 1\leq i<j\leq n,\ x_i>x_j\},\\
\coinv w&=\#\{(i,j) \mid 1\leq i<j\leq n,\ x_i<x_j\}.
\end{align*}
For a partition $\lambda$, let $\len(\lambda)$ denote its length, and write $\lambda^c$ for its conjugate partition.
\begin{Proposition}\label{prop:MacMahon}
There is a bijection between $\mathcal{W}_n$ and $\Sym_n \times \mathcal{P}_n$ sending 
$w \mapsto (\sigma,\lambda)$ such that

{\rm [P1]} $\max(w) = \des(\sigma) + \len(\lambda)$;

{\rm [P2]} $\sumt(w) = \maj(\sigma) + \sumt(\lambda)$;

{\rm [P3]} $\coinv(w) = \inv(\sigma)$.
\end{Proposition}
The above bijection is not new: it is essentially MacMahon's classical correspondence
\cite{MacMahon1913} (see also \cite{Carlitz1975, Foata2004H}). We reproduce its
construction here in order to keep the second proof self-contained.
Properties {\rm [P1]} and {\rm [P2]} are well known, since they are exactly what is
needed for the classical Carlitz-type identities. Property {\rm [P3]}, however,
which is the key to the $p$-refinement obtained in the present paper, does not
seem to have been recorded in the literature, to the best of our knowledge.

\begin{proof}
	We describe the bijection by means of an example. Let $n=8$ and 
	$w = 45412225 \in \mathcal{W}_n$.
Place $w$ in the first row and $\id = 12345678$ in the second row:
\begin{align*}
w  &=4 \ 5 \ 4 \ 1 \ 2 \ 2 \ 2 \ 5\\
\id &=1 \ 2 \ 3 \ 4 \ 5 \ 6 \ 7 \ 8
\end{align*}
We now repeatedly swap adjacent columns $\binom{a\ b}{x\ y}$ whenever $a<b$, and continue this process until the top row becomes a weakly decreasing word $\overline{w}$:
\begin{align*}
\overline{w}  &=5 \ 5 \ 4 \ 4 \ 2 \ 2 \ 2 \ 1\\
\sigma        &=2 \ 8 \ 1 \ 3 \ 5 \ 6 \ 7 \ 4\\
z             &=2 \ 2 \ 1 \ 1 \ 1 \ 1 \ 1 \ 0\\
\mu           &=3 \ 3 \ 3 \ 3 \ 1 \ 1 \ 1 \ 1
\end{align*}
At this stage, the bottom row is the permutation $\sigma$.
Next, define $z$ by
$z_j = \#\{i \mid j \le i < n,\ \sigma_i > \sigma_{i+1}\}$ for $j<n$, and $z_n = 0$.
Then set $\mu_j = \overline{w}_j - z_j$, and let $\lambda = \mu^c$ be the conjugate of the partition $\mu$.
In our running example, this gives $\lambda = 844$.

We check {\rm [P1]} and {\rm [P2]} using the facts that $\des(\sigma) = z_1$ and $\maj(\sigma) = \sumt(z)$.
For {\rm [P3]}, observe that at every intermediate stage, if the current top word is $w'$ and the bottom word is $\sigma'$, then
the quantity $\coinv(w') + \inv(\sigma')$ remains invariant under the column swaps. Thus
$\coinv(w) + \inv(\id) = \coinv(\overline{w}) + \inv(\sigma)$.
This implies {\rm [P3]} since $\inv(\id)=0$ and $\coinv(\overline{w})=0$.
\end{proof}
\begin{proof}[Second proof of Theorem \ref{th:util}]
By Proposition \ref{prop:MacMahon}, we obtain
\begin{align*}
\sum_{\sigma \in \Sym_n} t^{\des(\sigma)} q^{\maj(\sigma)} p^{\inv(\sigma)}
\times
	\sum_{\lambda \in \mathcal{P}_n} t^{\len(\lambda)} q^{\sumt(\lambda)} 
&=\sum_{w\in \mathcal{W}_n} t^{\max(w)} q^{\sumt(w)} p^{\coinv(w)}.
\end{align*}
Hence
\begin{align*}
A_n(t,q,p)  \times \frac{1}{(1-tq)\cdots (1-tq^n)} &=\sum_{\ell\geq 0}  t^{\ell}  \sum_{\substack{w\in \mathcal{W}_n\\ \max(w)=\ell}} q^{\sumt(w)} p^{\coinv(w)},
\end{align*}
or equivalently
\begin{align*}
\frac{A_n(t,q,p) }{(1-t)(1-tq)\cdots (1-tq^n)} &=\sum_{\ell\geq 0}  t^{\ell}  \sum_{\substack{w\in \mathcal{W}_n\\ \max(w)\leq\ell}} q^{\sumt(w)} p^{\coinv(w)}.
\end{align*}
Observe that, when $p=1$, the preceding identity reduces to Carlitz’s formula \eqref{eq:Carlitz}.
Extracting the coefficient of $t^\ell$ from both sides, we find
\begin{align*}
[t^{\ell}] \frac{A_n(t,q,p)}{(t;q)_{n+1}} 
&=\sum_{\substack{w\in \mathcal{W}_n\\ \max(w)\leq\ell}} q^{\sumt(w)} p^{\coinv(w)}\\
&=\sum_{\substack{\overline{w}\in \mathcal{D}_n\\ \max(\overline{w})\leq \ell}}q^{\sumt(\overline{w}) }
\sum_{w\in \Sym(\overline{w})}  p^{\coinv(w)}.
\end{align*}
Write $ \overline{w} = \ell^{n_{\ell}}\cdots 2^{n_2}1^{n_1}0^{n_0} $,
so that $n=n_0+n_1 +  \cdots + n_{\ell}$ and $\sumt(\overline w)=n_1 + 2 n_2 + \cdots + \ell n_{\ell}$.
By the classical combinatorial interpretation of the $q$-multinomial coefficient
\cite{Andrews1976, MacMahon1913}, the statistic $\coinv$ is distributed on the set
$\Sym(\overline w)$ of all rearrangements of $\overline w$ according to
$$
\sum_{w\in \Sym(\overline{w})}  p^{\coinv(w)}
= {n   \brack n_0,  \cdots, n_{\ell}}_p .
$$
Hence
\begin{equation*}
[t^{\ell}] \frac{A_n(t,q,p)}{(t;q)_{n+1}} 
= \sum_{n_0,  \cdots, n_{\ell} \vdash n }q^{n_1 + 2 n_2 + \cdots + \ell n_{\ell} }  {n   \brack n_0,  \cdots, n_{\ell}}_p.\qedhere
\end{equation*}
\end{proof}

We also recall the multinomial version of the $q$-Lucas theorem \cite{Adin2005GR_sigMahonians, Olive1965}.

\begin{Theorem}[Olive]\label{th:qLucas}
Let $n_0, \ldots, n_\ell$ be nonnegative integers and write
$n_0 = mk_0 + i_0, \ldots, n_\ell = mk_\ell + i_\ell$ with $0\leq i_0, \ldots, i_\ell \leq m-1$.
Then
\begin{equation*}
{n_0 + \cdots + n_\ell   \brack n_0, \cdots, n_{\ell}}_{\xi_m} 
=
\binom{ k_0 + \cdots + k_\ell  }{  k_0, \cdots, k_{\ell}}
	{  i_0 + \cdots + i_\ell \brack i_0, \cdots, i_{\ell}}_{\xi_m} .
\end{equation*}
In addition, if $i_0 + \cdots + i_\ell\geq m$,
then
$$
	{  i_0 + \cdots + i_\ell \brack i_0, \cdots, i_{\ell}}_{\xi_m}  =0.
$$
\end{Theorem}
Now we are ready to give the second proof of our main theorem, which is very short.
\begin{proof}[Second proof of Theorem \ref{th:main}]
Let $n = mk + i$ with $0 \le i \le m-1$.  
For nonnegative integers $n_0, n_1, \ldots, n_\ell$ satisfying $n_0 + \cdots + n_\ell = n$, write
$n_j = mk_j + i_j$ with $0 \le i_j \le m-1$.
Applying the first part of the $q$-Lucas theorem (Theorem \ref{th:qLucas}), we obtain
\begin{align*}
  &q^{n_1 + 2 n_2 + \cdots + \ell n_{\ell}} {n \brack n_0, \ldots, n_\ell}_{\xi_m} \\
  &\qquad =
  q^{m(k_1 + 2 k_2 + \cdots + \ell k_{\ell})} \binom{k}{k_0, \ldots, k_\ell}
  \times
  q^{i_1 + 2 i_2 + \cdots + \ell i_{\ell}} {i \brack i_0, \ldots, i_\ell}_{\xi_m}.
\end{align*}

We now sum this identity over all nonnegative integers $n_j$ $(0 \le j \le \ell)$ with $n_0 + n_1 + \cdots + n_\ell = n$.  
By the second part of Theorem \ref{th:qLucas},
$$
q^{n_1 + 2 n_2 + \cdots + \ell n_{\ell}} {n \brack n_0, \ldots, n_\ell}_{\xi_m} = 0
\quad \text{if } i_0 + i_1 + \cdots + i_\ell \ge m,
$$
so the sum effectively runs only over those sequences $(n_0,\ldots,n_\ell)$ for which $i_0 + \cdots + i_\ell \le m-1$.  
Equivalently, these sequences are characterized by the conditions $k_0 + \cdots + k_\ell = k$ and $i_0 + \cdots + i_\ell = i$.  
Thus,
\begin{align*}
  &\sum_{n_0, \ldots, n_\ell \vdash n} q^{n_1 + 2 n_2 + \cdots + \ell n_{\ell}} {n \brack n_0, \ldots, n_\ell}_{\xi_m} \\
  &\qquad =
  \sum_{k_0, \ldots, k_\ell \vdash k} q^{m(k_1 + 2 k_2 + \cdots + \ell k_{\ell})} \binom{k}{k_0, \ldots, k_\ell} \\
  &\qquad\quad \times
  \sum_{i_0, \ldots, i_\ell \vdash i} q^{i_1 + 2 i_2 + \cdots + \ell i_{\ell}} {i \brack i_0, \ldots, i_\ell}_{\xi_m}.
\end{align*}

Combining this relation with Theorem \ref{th:util} completes the proof of the main theorem.
\end{proof}

\section{Concluding remarks}\label{sec:remarks} 
We end the paper with three questions suggested by the referees, which we believe deserve
further investigation.

\begin{Question}
Wachs \cite{Wachs1992_involution} proved the even case \eqref{eq:Wachs} by means of an explicit
sign-reversing involution on permutations. Is there a direct combinatorial or bijective proof of
the missing odd case \eqref{eq:missing}, or of the order-$4$ case \eqref{eq:caseI}?
\end{Question}

\begin{Question}
Dong and Lin \cite[Theorem 1.11]{Dong2024L} obtained a refinement of the signed Eulerian
identities of Theorem \ref{th:DF} by the peak statistic. Is there a peak-statistic refinement of
the main theorem, Theorem \ref{th:main}?
\end{Question}

\begin{Question}
The signed Mahonian and signed Eulerian polynomials have been extensively extended to the
classical Weyl groups and, more generally, to the colored permutation groups $G(r,1,n)$
\cite{Adin2005GR_sigMahonians, Biagioli2006_sMah, Eu2021LL_sEulMah, Reiner1995_Des}. Can the
Hadamard product framework developed here be adapted to the cyclotomic Euler--Mahonian
polynomials of these groups?
\end{Question}

\medskip

{\it Acknowledgement}.\quad The author wishes to thank 
Kathy Q. Ji for carefully reading the manuscript and providing valuable suggestions and corrections that improved the presentation of this work. He is also grateful to the two anonymous referees for their careful reading and for many helpful remarks and suggestions.


\medskip

\bibliographystyle{plain}


\end{document}